\documentclass[a4paper,11pt]{amsart}
            \usepackage{amssymb,graphics}
             
          \newcommand{\vs}{\vspace{5pt}\smallbreak}
             \renewcommand{\,}{\kern.7pt}
                 \newcommand{\gY}{k}
                 \newcommand{\gN}{h}
               \newcommand{\wN}{\sigma}
                 \newcommand{\gM}{g}
               \newcommand{\wM}{\omega}
            \newcommand{\twM}{\tilde\wM_1}
                \newcommand{\kM}{g_0}
                 \renewcommand{\k}{r}
                 \renewcommand{\l}{s}
                  \newcommand{\s}{u}
                  \newcommand{\x}{z}
                 \renewcommand{\v}{y}
                 \renewcommand{\u}{x}
              \newtheorem{lemma}{Lemma}
            \newtheorem{prop}{Proposition}
              \newtheorem{thm}{Theorem}
             \newtheorem{cor}{Corollary}
              \theoremstyle{definition}
               \newtheorem{rem}{Remark}
            
               \newtheorem{ex}{Example}

             \newcommand{\R}{\mathbb{R}}
             \newcommand{\C}{\mathbb{C}}
             \newcommand{\Z}{\mathbb{Z}}
             \newcommand{\T}{\mathbb{T}}
            \newcommand{\cL}{{\mathcal L}}
            \newcommand{\cM}{{\mathcal M}}
                 \newcommand{\X}{X}
                \newcommand{\ab}{(a,b)}
            \newcommand{\Ai}{\mathrm{Ai}}
               \newcommand{\we}{\wedge}
               \newcommand{\na}{\nabla}
               \newcommand{\al}{\alpha}
               \newcommand{\la}{\lambda}
               \newcommand{\Om}{\Omega}
            \newcommand{\OO}{\Om\we\ol\Om}
               \newcommand{\ga}{\gamma}
              \newcommand{\pd}{\partial}
              \newcommand{\La}{\Lambda}
          \renewcommand{\Re}{\mathop{\mathfrak{Re}}}
          \renewcommand{\Im}{\mathop{\mathfrak{Im}}}
        \newcommand{\q}{\quad}\newcommand{\qq}{\qquad}
          \newcommand{\y}{\\[4pt]}\newcommand{\yy}{\\[6pt]}
           \newcommand{\ft}[2]{\hbox{$\frac#1#2$}}
            \newcommand{\fp}{\frac1{2\pi}}
               \newcommand{\ds}{\displaystyle}
       \newcommand{\alt}{\raise1.5pt\hbox{$\bigwedge$}}
          \newcommand{\cd}{\>\hbox{\Large$\cdot$}\>}
                \newcommand{\+}{\!+\!}
               \renewcommand{\=}{\!=\!}
               \newcommand{\hk}{hyperk\"ahler }
             \newcommand{\Ka}{K\"ahler }
         \newcommand{\acs}{almost-complex structure }
            \newcommand{\aKa}{almost-\Ka}
            \newcommand{\eps}{\varepsilon}
          
               \newcommand{\ba}{\begin{array}}
            \newcommand{\ea}{\end{array}}
             \newcommand{\ol}{\overline}
               \newcommand{\ti}{\times}
              \newcommand{\ot}{\otimes}
   \def\FFF{\vs\vs\centerline{\scalebox{.75}{\includegraphics{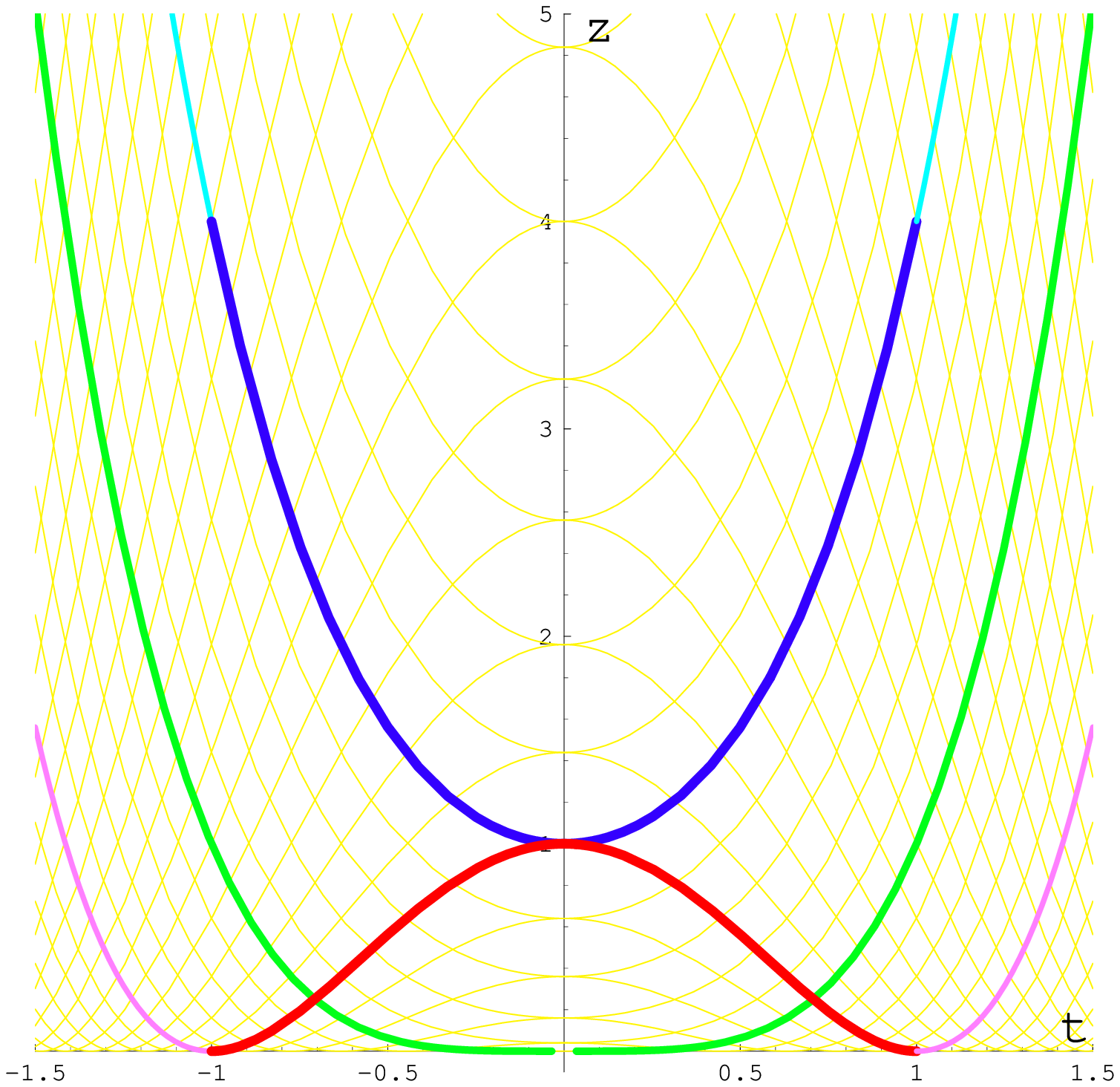}}\vs}
            \par\centerline{Figure}\vs\vs}
                   
\begin{document}

\title{K\"ahler reduction of metrics with holonomy G$_{\mathbf2}$}

\author[V. Apostolov]{Vestislav Apostolov}
\address{Vestislav Apostolov \\ D{\'e}partement de Math{\'e}matiques\\ UQAM\\
C.P.~8888 \\ Succ.\ Centre-ville \\ Montr{\'e}al (Qu{\'e}bec) \\ H3C 3P8 \\
Canada}
\email{apostolo@math.uqam.ca}
\author[S.~Salamon]{Simon Salamon}
\address{Simon Salamon \\ Dipartimento di Matematica\\ Politecnico di
Torino\\Corso Duca degli Abruzzi 24\\10129 Torino, Italy}
\email{salamon@calvino.polito.it}

\maketitle

\setcounter{section}0
\section*{Introduction}

There are now many explicitly known examples of metrics with holonomy group
equal to $G_2$, the simplest of which admit an isometry group with orbits of
codimension one. A metric with holonomy $G_2$ on a smooth 7-manifold $Y$ is
characterized by a 3-form and a 4-form that are interrelated and both
closed. If this structure is preserved by a circle group $S^1$ acting on $Y$
(which cannot then be compact), the quotient $Y/S^1$ has a natural symplectic
structure. A geometrical description of such quotients was carried out by
Atiyah and Witten \cite{AtWMTD} when $Y$ is one of three original manifolds
with a complete metric of holonomy $G_2$ described in \cite{BrSCCM}. The
symplectic structure is also defined on the image of the fixed points, which
in each of these cases is a Lagrangian submanifold $L$. The embedding of $L$
in a neighborhood of $Y/S^1$ is believed to approximate the geometry of a
special Lagrangian submanifold of a Calabi-Yau manifold, and this is
consistent with models of special Lagrangian submanifolds of $\R^6$ described
for example in \cite{JoyUIS}.

This work motivates a general investigation of the quotient $N=Y/S^1$ of a
manifold with holonomy $G_2$ by an $S^1$ action. Initial results in this
direction can be found in \cite{ChSITS,CGLASH}, and in this paper we pursue the
theory under the simplifying assumption that the $S^1$ action is free. It is an
elementary fact that $N$ has, in addition to its symplectic 2-form, a natural
reduction to $SU(3)$. Whilst this structure cannot be torsion-free (i.e.\
Calabi-Yau) if $Y$ is irreducible, there are non-trivial examples in which the
associated \acs is integrable, so that $N$ is K\"ahler.  This paper is devoted
to an investigation of such a situation, which turns out to be surprisingly
rich. Our study began with the realization that in this case, a parameter
measuring the size of the fibers of $Y$ generates both a Killing and a Ricci
potential for $N$ (Proposition~1 and Corollary~2), exhibiting a link with the
theory of the so-called Hamiltonian 2-forms \cite{ACGHFK}.

We first describe the induced $SU(3)$ structure on $N$ in \S1, and then pursue
consequences of the integrability condition. In particular, we prove that when
$N$ is K\"ahler, an infinitesimal isometry $U$ of the $SU(3)$ structure is
inherently defined. The situation is reminiscent of the study of
Einstein-Hermitian 4-manifolds in which Killing vector fields appear
automatically \cite{AG0,DE,LeBEMC}. We explain in \S2 that $U$ can be used to
obtain a \Ka quotient of $N$, consisting of a 4-manifold equipped with a
1-parameter family of smooth functions and 2-forms, satisfying a coupled
second-order evolution equation (Theorem~\ref{thm1}). The procedure can be
reversed so as to construct metrics with holonomy (generically equal to) $G_2$
from a 4-manifold $M$ with the appropriate structure.

A fundamental construction of a holonomy $G_2$ metric starting from a $\T^2$
bundle over a \hk 4-manifold was discovered by Gibbons, L\"u, Pope and Stelle
\cite{GLPSDW}, and towards the end of \S3 we exhibit our reverse procedure as a
generalization of the above. It can be further improved by introducing an
anti-self-dual 2-form with the opposite orientation to the \Ka one
(Theorem~\ref{thm2}). This leads to the construction of new examples of $G_2$
metrics based on the examples of Ricci-flat \aKa metrics
\cite{ACGGWS,ArmAAK,NuPFDE}. We point out that there also exist constructions
of holonomy $G_2$ metrics from higher dimensional \hk manifolds (see for
example \cite{AcWCFM}).

When $M$ is $\T^4$, the basic examples are modelled on nilpotent Lie groups and
fall into three types, special cases of which were also mentioned in
\cite{ChSITS}. We explain how to show that these are irreducible and therefore
have holonomy equal to $G_2$. In \S4, we analyse the holonomy of these metrics
by restricting them to hypersurfaces so as to induce an $SU(3)$ structure that
evolves according to equations studied by Hitchin \cite{HitSFS}. Our examples
provide perhaps the simplest instance of this phenomenon other than the case of
nearly-\Ka manifolds. The different $SU(3)$ structures obtained in this way, an
integrable one on a quotient and non-integrable ones on hypersurfaces, are
naturally linked via the common 7-manifold.

In the final section, we study a more general class of solutions of our system,
formulated in terms of the complex Monge-Amp\`ere equation. This leads to both
an abstract existence theorem, some special solutions, and an explicit final
example.

Our assumption that the symplectic manifold $Y/S^1$ be \Ka leads to a (local)
action on $Y$ not just by $S^1$ but by the torus $\T^2$
(Corollary~\ref{newcor}), though examples in \S4 exhibit $\T^2$ actions for
which the \Ka assumption fails (see also \cite{KDLI6B,SanCG2}). We hope
nonetheless that similar methods will lead to the classification of metrics
with holonomy $G_2$ admitting a $\T^2$ action.\vs

\noindent\textbf{Acknowledgments.} This material was first developed at the
joint AMS--UMI conference in Pisa in 2002. The first author was supported in
part by NSERC grant OGP0023879 and by NSF grant INT-9903302. The second author
is a member of EDGE, Research Training Network HPRN-CT-2000-00101, supported by
the European Human Potential Programme, and a preliminary version of the
results was presented at its mid-term conference. The authors are grateful to
the referee for useful remarks, and to D.~Conti for subsequent checking of some
curvature computations.

\vs\section{The first reduction}

In this section we study general properties of a 7-dimensional Riemannian
manifold $(Y,\gY)$ endowed with a torsion-free $G_2$ structure and an
infinitesimal isometry $V$. The $G_2$ structure on $Y$ is defined by an
admissible 3-form $\varphi$, which itself determines the Riemannian metric
$\gY$, an orientation and Hodge operator $*$. The torsion-free condition is
equivalent to the closure of both $\varphi$ and the 4-form $*\varphi$, and
amounts to asserting that the holonomy group of $k$ is contained in $G_2$.
This theory is elaborated in the standard references \cite{BryMEH,FeGRMS,
JoyCMS,SalRGH}

We denote by $(N,\check\gY)$ the Riemannian quotient of $(Y,\gY)$, so that $N$
is a 6-dimensional manifold formed from the orbits of the Killing vector field
$V$. The considerations are local and hold in a suitable neighborhood of any
point of $Y$ where the vector field $V$ is non-zero. The Riemannian metric
$\check\gY$ is induced by the formula \[\gY=\check\gY+t^{-2}\eta\ot\eta,\]
where $t=\|V\|^{-1}_\gY= \gY(V,V)^{-1/2}$ and $\eta=\gY(\cd,t^2V)$ is the
1-form dual to $t^2V$ so that $\eta(V)=1$. In the case when $Y$ can be realized
as a principal $S^1$ bundle over $N$ (the fibers being the closed orbits of
$V$), $\eta$ is nothing but a connection 1-form.

Since $V$ is a Killing vector field, $dt(V)=0$ and $t$ descends to a function
on $N$. The 2-form \[\wN=\imath_V\varphi\] is horizontal in the sense that
$\imath_V\wN=0$, and since $V$ preserves the $G_2$ structure we have \[\ba{c}
d\wN=d(\imath_V\varphi)=\cL_V\varphi=0,\y \cL_V\wN = d(\imath_V\wN)=0.\ea\]
Thus, $\wN$ is closed and $V$-invariant, and therefore the pullback of a
closed 2-form on $N$, again denoted by $\wN$. (We identify functions and forms
on $N$ with their pullbacks to $Y$ throughout, and $\cL$ and $\imath$ stand
for the Lie derivative and the interior product.) Since $\wN$ is
non-degenerate transverse to the fibers of $Y$, it defines a symplectic form
on $N$.

We now choose to express the forms characterizing the $G_2$ structure on $Y$
as \begin{equation}\label{express}\ba{rcl}\varphi &=&\wN\we\eta\ +\ t\,\Psi^+
\yy*\varphi &=& \Psi^-\we\eta+\ft12t^2\wN\we\wN,\ea\end{equation} where
\[\Psi^-= -\imath_V (*\varphi).\] If $\eta$ were replaced by the unit form
$t^{-1}\eta$ then both $\wN$ and $\Psi^-$ would need to be multiplied by a
compensatory factor of $t$ if the left-hand sides of (\ref{express}) are to
remain the same. This explains why $\wN\we\wN$ appears in (\ref{express}) with
a coefficient of (one half) $t^2$, and $\Psi^+$ appears with a coefficient of
$t$ so as to have the same norm as $\Psi^-$.  The rescaled Riemannian metric
\begin{equation}\label{gNgY}\gN=t^{-1}\check\gY\end{equation} is compatible in
the sense that the skew-symmetric endomorphism $J$ defined by
\begin{equation}\label{gJ}\gN(J\cd,\cd)=\wN(\cd,\cd)\end{equation} is an \acs
on $N$. The triple $(\gN,\wN,J)$ then defines an \emph{\aKa}structure on $N$,
though the qualification `almost' can be deleted when $J$ is integrable.

Just as for $\wN$, the $3$-form $\Psi^-$ is closed and basic (meaning
$\imath_V\Psi^- =0$, $\cL_V\Psi^-=0$). Thus, it too is the pullback of a real
form on $N$. This has type $(3,0)\+(0,3)$ with respect to $J$, so that
$\Psi^-(J\cd,J\cd,\cd)=-\Psi^-(\cd,\cd,\cd)$. The theory of $G_2$ structures
then implies that \begin{equation}\label{Psi+}\Psi^+(\cd,\cd,\cd)=
\Psi^-(J\cd,\cd,\cd),\end{equation} and the two real 3-forms combine to define
a complex form \begin{equation}\label{Psi} \Psi=\Psi^++ i\Psi^-\end{equation}
of type $(3,0)$ with respect to $J$.

Unlike $\Psi^-$, the 3-form $\Psi^+$ is not in general closed; indeed,
$d\Psi^+$ can be identified with the Nijenhius tensor of $J$, the obstruction
to the integrability of $J$:\vs

\begin{lemma}\label{lem1} $\Psi^+$ is closed if and only if the \acs $J$ is
integrable.\end{lemma}

\begin{proof} Since $\Psi^-$ is already closed, the exterior derivative
$d\Psi^+$ of (\ref{Psi}) is real. But if $J$ is integrable so that $(N,J)$ is
a complex manifold, this real 4-form has type $(3,1)$ and therefore vanishes.
Conversely, if $d\Psi^+=d\Psi^-=0$ then writing (\ref{Psi}) locally as a wedge
product of $(1,0)$-forms $\al^i$ shows that $(d\al^i)^{0,2}=0$ for $i=1,2,3$,
and $J$ is integrable.\end{proof}\vs

The \aKa structure corresponds to a reduction to $U(3)$ at each point of
$N$. Specification of the non-zero element (\ref{Psi}) in the space
$\La^{3,0}$ is precisely what is needed to reduce the structure to $SU(3)$,
but one should first rescale given that $\|\Psi^\pm\|_\gN$ is not in general
constant. The $3$-forms \[\psi^\pm= t^{-1/2}\Psi^\pm\] have norm equal to
2. They are subject to the compatibility equations
\begin{equation}\label{comp1}\wN\we\psi^\pm=0, \end{equation} and
\begin{equation}\label{comp2} \psi^+\we\psi^-=\ft23\wN^3=4\mathrm{vol}_\gN,
\end{equation} consistent with (\ref{gJ}), and give rise to an
$SU(3)$ structure with underlying metric $\gN$. Notation for the 3-forms
is now consistent with that of \cite{ChSITS}. Since
$\psi^++i\psi^-\in\La^{3,0}$, we have \[(\ga
+iJ\ga)\we(\psi^++i\psi^-)=0\] for any 1-form $\ga$, where $J$ acts on
1-forms by $J\ga(\cd)=-\ga(J\cd)$. Whence \begin{equation}
\label{Jpsi}\ga\we\psi^+=J\ga\we\psi^-.\end{equation} This last equation will
be useful for calculations.

The real 3-forms $\psi^+,\psi^-$ have a common stabilizer group $SL(3,\C)$ at
each point of $N$, and either one determines the \acs $J$ \cite{HitGTF}. From
this point of view, the further reduction to $SU(3)$ is achieved by a 2-form
$\wN\in\La^{1,1}$ such that $\wN(\cd,J\cd)=\gN(\cd,\cd)$ is positive definite.
The $SU(3)$ structure is therefore fully determined by (say) $\psi^+$ and
$\wN$; in the new notation Lemma~\ref{lem1} reads as\vs

\begin{cor}\label{cor1} The \acs $J$ is integrable if and only if \[\na\psi^+
=-\ft12(d^c\log t)\ot\psi^-,\q\na\psi^-=\ft12 (d^c\log t)\ot \psi^+,\] where
$\na$ denotes the Levi-Civita connection of $\gN$ and $d^c f= J\,df$.\end{cor}

\begin{proof} Given the above equations, \[d\Psi^+=d(t^{1/2}\psi^+)=\ft12
t^{1/2}(d\log t\we\psi^+-d^c\log t\we\psi^-),\] but the right-hand side is zero
by (\ref{Jpsi}). Thus, $J$ is integrable by Lemma \ref{lem1}.

Conversely, if $J$ is integrable then $(\gN,J)$ is \Ka and $\na_U\psi^\pm$ is a
form of type $(3,0)\+(0,3)$ for any vector field $U$. Since $\psi^+$ and
$\psi^-$ are mutually orthogonal and of constant norm $2$, \[\na\psi^+=\ga\ot
\psi^-,\q\na \psi^-=-\ga\ot\psi^+\] for some 1-form $\ga$. Using (\ref{Jpsi})
again, we obtain \[\ba{rcl} d\Psi^+ &=& d(t^{1/2}\psi^+)=\ft12t^{1/2}(d\log t
\we\psi^++2\ga\we\psi^-)\y &=& \ft12t^{1/2}(-2J\ga+d\log t)\we\psi^+.\ea\] The
claim follows by Lemma~\ref{lem1} and the injectivity of the linear map
$\La^1\to\La^4$ defined at each point by wedging with $\psi^+$.\end{proof}\vs

\begin{cor}\label{cor2} If $J$ is integrable, the Ricci form $\kappa$ of
the \Ka manifold $(N,\gN,J)$ is given by \[\kappa=\ft12dd^c\log t.\]\end{cor}

\begin{proof} This is again an immediate consequence of the integrability
criterion, which implies that the $(3,0)$-form $\Psi=\Psi^++i\Psi^-$ is
closed. It follows that $\Psi$ is a holomorphic section of the canonical bundle
$\La^{3,0}$ and the Ricci form is given by \[\kappa=i\pd\ol\pd\log(
\|\Psi\|^2_\gN)=\ft12dd^c\log(\|\Psi\|^2_\gN),\] and the result follows from
(\ref{comp2}).\end{proof}\vs

In the case when $J$ is integrable, it is now possible to formulate a complete
system of conditions for the $SU(3)$ structure on $N$ to arise from the
quotient of a torsion-free $G_2$ structure.\vs

\begin{prop}\label{criterion} Let $(N,\gN,J,\wN)$ be a \Ka manifold of real
dimension 6, endowed with a compatible $SU(3)$ structure defined by a
$(3,0)$-form $\psi^++i\psi^-$. Then this structure is obtained as the quotient
of a 7-dimensional manifold $Y$ with a torsion-free $G_2$ structure by a
nontrivial infinitesimal isometry $V$ if and only if \begin{enumerate}
\item $d\psi^+= -\frac12(d^c\log t)\we\psi^-$ and $d\psi^-=\frac12(d^c\log t)
\we\psi^+$ for a smooth positive function $t$, and \item the Hamiltonian vector
field $U$ on $(N,\wN)$, corresponding to $-t$, is an infinitesimal isometry of
the $SU(3)$ structure.\end{enumerate} In this case, the corresponding
7-manifold $Y$ is locally $\R\ti N$ and the metric $\gY$, infinitesimal
isometry $V$ and $G_2$ invariant forms $\varphi,\, *\varphi$ are given by
\begin{equation}\label{phi*phi}\ba{c}\ds\gY=t\gN+t^{-2}\eta \ot\eta,\qq
V=\frac\pd{\pd \v}\yy \varphi=\wN\we\eta+t^{3/2}\psi^+,\qq *\varphi=t^{1/2}
\psi^-\we\eta +\ft12t^2\wN\we\wN,\ea\end{equation} where $\v$ is a variable for
the $\R$ factor, and $\eta=d\v+\eta_N$ is a 1-form on $\R\ti N$ for which
\begin{equation}\label{theta}d\eta_N=-t^{1/2}(\imath_U\psi^+).
\end{equation}\end{prop}

\begin{proof} Condition (i) is necessary by Corollary \ref{cor1}, and the
equalities (\ref{phi*phi}) reflect the earlier definitions of
$\gN,\,\wN,\, \psi^\pm$. Moreover, $\pd/\pd {\v}$ is identified
with the Killing field $V$ (so that $\gY(V,V)=t^{-2}$) and $\eta$
is the corresponding connection form. Using (i) and (\ref{Jpsi})
in that order, we obtain \[\ba{rcl}0\ =\ d\varphi &=&\wN\we
d\eta+\ft32t^{1/2}dt\we\psi^+-\ft12t^{1/2}Jdt\we\psi^-\y &=&\wN\we
d\eta+t^{1/2}dt\we\psi^+.\ea\] Consequently, $(d\eta)^{1,1}=0$ and
(\ref{theta}) follows from the fact that $\imath_U\wN=-dt$.

We now have \[0=-d(d\eta)=d(t^{1/2}\imath_U\psi^+)=d(\imath_U\Psi^+)=\cL_U
\Psi^+=t^{1/2}\cL_U\psi^+.\] Since $U$ is Hamiltonian, we also have \[\cL_U
\wN =0,\] i.e.\ $U$ is an infinitesimal isometry for the pair $(\wN,\psi^+)$,
and therefore for the $SU(3)$ structure (using \cite{HitGTF}). Reversing the
above arguments, one can check directly that (\ref{phi*phi}) and (\ref{theta})
define a torsion-free $G_2$ structure on $Y=\R\ti N$.\end{proof}\vs

\begin{cor}\label{newcor} Under the hypothesis of
Proposition~\ref{criterion}, the horizontal lift of $U$ to $(Y,k)$
is an infinitesimal isometry of the $G_2$ structure $\varphi$,
which commutes with $V$.
\end{cor}
\begin{proof} This is an immediate consequence from
Proposition~\ref{criterion}.
\end{proof}

\begin{rem} {\rm It follows from Proposition \ref{criterion} that when $t$ is
constant, $(Y,\gY)$ is locally the Riemannian product of a Calabi-Yau
$6$-manifold with $\R$. In this case the holonomy group of $\gN$ lies in
$SU(3)$. In general, the failure of the holonomy to reduce to $SU(3)$ is
measured by a torsion tensor \[ \tau\in T^*_nN\ot\frac{\mathfrak{so}(6)}
{\mathfrak{su}(3)}\] determined by $d\wN,d\psi^+,d\psi^-$. Whereas $\tau$ has
a total of 42 components, exactly two thirds of these will always vanish on
$N$ and $\tau$ is determined by the remaining 14 tracefree components of
$d\psi^+$, or 6 in the \Ka case \cite{ChSITS}.}\end{rem}

\vs\section{A second reduction}

Proposition~\ref{criterion} is the key ingredient for performing a further
quotient via the infinitesimal isometry $U$. To carry this out, we shall assume
from now on that $t$ is not constant and that $J$ is integrable.  Thus, $U$ is
a non-trivial infinitesimal isometry of the $SU(3)$ structure on $N$ determined
by the pair $(\gN,\psi^+)$.

Denote by $(M,J_1)$ the `stable' or holomorphic quotient of $(N,J)$, defined
at least locally as the complex two-dimensional manifold of holomorphic leaves
of the foliation generated by $\Xi= U-iJU$. Here we exploit the fact that
$\Xi$ is a holomorphic vector field on $(N,J)$. Since $J$ is integrable,
$\Psi=\Psi^++i \Psi^-$ is a closed $(3,0)$-form, and we set \[\Om=\ft12
\imath_\Xi\Psi.\] Since $\imath_\Xi\Om=0$, it follows from (\ref{Psi+}) that
$\Om$ is closed and the pullback of a holomorphic symplectic form on $M$,
again denoted by $\Om$. The real closed 2-forms \[\wM_2=\Re\Om,\q
\wM_3=\Im\Om\] on $M$ (that pull back to $\imath_U\Psi^+$, $\imath_U\Psi^-$
respectively on $N$) satisfy \begin{equation}\label{23}
\wM_2\we\wM_2=\wM_3\we\wM_3,\qq\wM_2\we\wM_3 =0.\end{equation} The complex
structure $J_1$ can now be determined by the formula
\[\wM_2(\cd,\cd)=\wM_3(J_1\cd,\cd),\] and equation (\ref{theta}) in
Proposition \ref{criterion} reads \begin{equation}\label{deta}d\eta=-\wM_2.
\end{equation}\medbreak

Set $\s=\|U\|_\gN^{-2}$ and let $\xi$ be the 1-form $\gN$-dual to
the vector field $\s U$, so that $\imath_U\xi=1$ and $\xi=\s Jdt$.
The 3-forms $\Psi^\pm$ are completely determined by (\ref{theta})
which forces them to be the real and imaginary components of
$(\wM_2+i\wM_3)\we(\xi+iJ\xi)$. Thus
\begin{equation}\label{final1}\ba{c}
\psi^+=t^{-1/2}(\wM_2\we\xi+\s\,\wM_3\we dt),\yy
\psi^-=t^{-1/2}(\wM_3\we\xi-\s\,\wM_2\we dt).\ea\end{equation} They are
$U$-invariant in accordance with (ii) in Proposition~\ref{criterion}.

For any regular value of the momentum map $-t$, the stable quotient $(M,J_1)$
of $(N,J)$ can be identified with the \emph{symplectic} quotient $(M,\twM(t))$
of $(N,\wN)$ generated by the vector field $U$. In this way, we obtain the
\emph{\Ka} quotient $(M,\gM(t),\twM(t),J_1)$. In this correspondence,
\begin{equation}\label{final2}\ba{c} \wN=\twM(t)+dt\we\xi,\yy\ds\gN=\gM(t)+
\s^{-1}\xi\ot\xi+\s\,dt\ot dt,\ea\end{equation} so that the equation
(\ref{comp2}) reduces to \begin{equation}\label{Mcompatibility}
t\,\twM(t)\we\twM(t)=\ft12\s\,\,\OO=\s\,\wM_2\we\wM_2=\s\,\wM_3\we\wM_3.
\end{equation} To ease the notation, we shall below omit the explicit
dependence of $\twM=\twM(t)$ on $t$ except on occasions for emphasis.

We now denote by $P$ the (locally defined) space of orbits of $U$, so that $N$
can be thought as an $\R$ bundle over $P$ with connection 1-form $\xi$.
Locally, $N=\R\ti P$, and introducing a variable $\u$ for the $\R$ factor, we
may write \[\xi=d\u+\xi_P,\qq U=\frac\pd{\pd \u}\] for some 1-form $\xi_P$ on
$P$. The space of orbits of the vector field $JU$ on $P$ is the stable quotient
of $(N,J)$, whereas the symplectic quotients of $(N,\wN)$ are identified with
the level sets of $t$ in $P$. Using a local description $P=\R^+\ti M$ in which
the $\R^+$-factor corresponds to $t$, we may regard $\s$ as a function on $M$
for each value of $t$. In these terms, we have \[\ba{c}d\xi=d_P\xi_P=\al_M\we
dt+\beta_M,\yy d_P\s=\s'dt+d_M\s,\yy d\twM=\twM'\we dt,\ea\] where
$\al_M=\al_M(t),\beta_M=\beta_M(t)$ are 1-parameter families of forms on $M$,
$\ '$ denotes $\pd/\pd t$ and $d,d_P,d_M$ denote exterior derivative on
$N,P,M$.

Differentiating the first relation gives \[d_M\beta_M=0,\q\beta_M'=-d_M\al_M.
\] Using the formula (\ref{final2}) for the symplectic form $\wN$ yields \[
\twM'=\beta_M,\] and it follows that $\beta_M$ has type $(1,1)$ relative to
$J_1$. Using (\ref{final1}) gives \[\ba{rcl} d\psi^+ &=&-\ft12t^{-3/2}dt\we
\wM_2\we\xi+t^{-1/2}\wM_2\we d\xi+ t^{-1/2}d_M\s\we\wM_3\we dt\y &=&
-\ft12t^{-1}Jdt\we\psi^-+t^{-1/2}(\wM_2\we\al _M+d_M\s\we\wM_3)\we dt,\ea\]
since $\wM_2\we\beta_M=0$. In view of Proposition~\ref{criterion}(i),
$\wM_2\we\al_M+d_M\s\we\wM_3=0$, whence \[\al_M=J_1d\s=d^c_M\s,\] and
everything can be expressed in terms of $\twM$ and $\s$.

In summary,\vs

\begin{thm}\label{thm1} Let $(Y,\varphi)$ be a 7-manifold with a torsion-free
$G_2$ structure, admitting an infinitesimal isometry $V$. Suppose that the norm
of $V$ is not constant, and let $y\in Y$ be a point where $V$ does not vanish.
Let $\mathcal V$ be a neighborhood of $y$, such that the space of orbits of $V$
in $\mathcal V$ is a manifold $N$ and suppose that the \aKa structure on $N$ is
in fact K\"ahler. Then, there exists a 4-dimensional manifold $M$ endowed with
a complex structure $J_1$, a complex symplectic form $\Om=\wM_2+i\wM_3$, and
1-parameter families of \Ka 2-forms $\twM=\twM(t)$ and positive functions
$\s=\s(t)$ on $M$, satisfying the relations \begin{equation}\label{flow}
\twM''=-d_Md^c_M\,\s\end{equation}\vspace{-15pt}\begin{equation}
\label{Tcompatible}t\,\twM\we\twM=\ft12\s\,\,\OO.\end{equation}

On a sufficiently small neighborhood of $y$, $(\varphi,V)$ is equivariantly
isometric to the torsion-free $G_2$ structure \begin{equation}\label{phi}
\ba{rcl}\varphi &=& \twM\we(d\v+\eta_N)+dt\we(d\u+\xi_P)\we(d\v+\eta_N)\y
&&\hspace{100pt}+\,t\Big(\wM_2\we(d\u+\xi_P)+\s\,\wM_3\we dt\Big),\ea
\end{equation} on $\R_t^+\ti\R^2_{\u,\v}\ti M$ endowed with the
infinitesimal isometry $V =\pd/\pd \v$, where \begin{enumerate}
\item[] $d_M$ denotes the differential on $M$, and $d^c_M=J_1\circ d_M$;
\item[] $t>0$ is the variable on the $\R_t^+$-factor; \item[] $(\u,\v)$ are
standard coordinates on $\R^2_{\u,\v}=\R_{\u}\ti\R_{\v}$;
\item[] $\eta_N$ is a 1-form on $N=\R_t^+\ti\R_{\u}\ti M$ with
$d\eta_N=-\wM_2$;
\item[] $\xi_P$ is a 1-form on $P=\R_t^+\ti M$ with $d\xi_P=(d^c_M\s) \we
dt+\twM'$.\end{enumerate}\end{thm}\vs

\begin{rem}\label{irreducible} By redefining the local coordinate $x$, one 
can assume (without loss) that $\eta_N$ is in fact a 1-form on $M$. In this
case the $G_2$ structure $\varphi$ is invariant under $\frac{\partial}{\partial
x}$ as well, and $\frac{\partial}{\partial x}$ is identified with the Killing
vector field defined in Corollary~\ref{newcor}.

It is not difficult to see that for {\it generic} data on $M$, the holonomy
group of the $G_2$ structure (\ref{phi}) is {\it equal} to $G_2$. Indeed, the
general theory of holonomy groups (see e.g.~\cite{SalRGH}) implies that if the
holomomy group of $(Y,k)$ were strictly less than $G_2$, then there would exist
a non-trivial parallel vector field $\X$ on $(Y,k)$ commuting with
$V=\frac{\partial}{\partial y}$; it would therefore come from an infinitesimal
isometry of the $SU(3)$ structure on $N$ (still denoted by $\X$), which
preserves the level sets of $t$ and commutes with $U=\frac{\partial}{\partial
x}$. The equations for $d\eta_N$ and $d\xi_P$ imply that there is no parallel
vector field in span$\{ \frac{\partial}{\partial x},\frac{\partial}{\partial y}
\}$ unless $d_Mu=0$ and ${\tilde \omega}' =0$ (a situation which we shall
exclude below), thus showing that $\X$ is in general different from
$U$. Therefore, $\X$ must come from a (real) holomorphic vector field on
$(M,J)$ which preserves the K\"ahler metrics $\omega(t)$ for each $t$. In
particular, if we assume that $(M,J_1,\Om,\twM(t),\s(t))$ does not admit any
infinitesimal isometry and either $d_Mu\ne0$ or ${\tilde \omega}'\ne0$, then
we know that the holonomy group of $(Y,\varphi)$ must equal $G_2$.\end{rem}

It is important to note that the triple $(\twM(t),\wM_2,\wM_3)$ appearing in
Theorem~\ref{thm1} does not in general constitute a \hk structure. Indeed, by
(\ref{Tcompatible}), the holomorphic section $\Om$ of the canonical bundle
$\La^{2,0}M$ satisfies \[\|\Om\|^2_{\gM}=t\s^{-1},\] and the Ricci form
$\kappa$ of the \Ka metric $\twM(t)$ is given by
\begin{equation}\label{ricM}\kappa=\ft12d_Md_M^c (\log
t-\log\s)=-\ft12d_Md_M^c\log\s.\end{equation} In the next section, we shall
however explain how the simplest case does correspond to a \hk situation.

\vs\section{Constant solutions}

A careful inspection of the proof of Theorem~\ref{thm1} shows that the process
can be inverted so as to construct a torsion-free $G_2$ structure from a
4-manifold $M$ with a complex symplectic structure $(J_1,\Om)$ together with a
1-parameter family $(\twM(t),\s(t))$ of \Ka forms and smooth functions
satisfying (\ref{flow}) and (\ref{Tcompatible}). In this section, we shall
carry out this inverse construction explicitly in the case in which $\s$ really
is just a function of $t$, so that $d_M\s=0$.

The above assumption reduces (\ref{flow}) to \[\twM=\wM +t\wM'\] for some
closed (1,1)-forms $\wM,\wM'$ on $M$. Consider the real symmetric bilinear form
$B$ defined by \begin{equation}\label{B}\al\we\beta=\ft12B(\al,\beta)\OO
\end{equation} on the space of 2-forms. Restricting $B$ to the subspace
$\langle\wM,\wM'\rangle$ and diagonalizing, we may write \begin{equation}
\label{example1} \twM=(p+qt)\wM_0+(\k+\l t)\wM_1,\end{equation} where
\begin{equation}\label{AK-HK}\wM_0\we\wM_0=-\ft12\eps\,\OO,\q\wM_1\we\wM_1=
\ft12\OO,\q\wM_0\we\wM_1 =0,\end{equation} $\eps$ is $0$ or $1$, and $p,q,
\k,\l$ are constants satisfying\begin{equation}\label{inequality}
\k+\l t>|p+qt|\end{equation} to ensure the overall positivity of $\twM$.
>From (\ref{example1}) and (\ref{Tcompatible}) we have \begin{equation}
\s=t\Big( (\k+\l t)^2-(p+qt)^2\Big).\end{equation} Note that $\eps=0$ in
(\ref{AK-HK}) if and only if $\wM_0=0$; to avoid redundancy we declare that
$p=q=0$ in this case.

The real symplectic forms $\wM_1,\wM_2,\wM_3$ satisfy the usual compatibility
relations \begin{equation}\label{HK}\ba{c}\wM_i\we\wM_i=\wM_j\we\wM_j,\y
\wM_i\we\wM_j=0,\q i\ne j,\ea\end{equation} extending (\ref{23}). It is
well known (\cite{HitSDE,SalRGH}) that they then determine a \hk structure,
consisting of \begin{enumerate} \item complex structures $J_1,J_2,J_3$
satisfying $J_i\circ J_j+J_j\circ J_i= -2\delta_{ij}{\rm Id}$; \item a
Riemannian metric $\kM$ and associated Levi-Civita connection relative to
which the $J_i$ are all orthogonal and parallel. \end{enumerate}

When $p^2+q^2>0$, in addition to the \hk structure, the $(1,1)$ form $\wM_0$
defines an \acs $I$ on $M$, such that $\wM_0(\cd,\cd)=\gM_0(I\cd,\cd)$ as in
(\ref{gJ}). It follows that $(\gM_0,I,\wM_0)$ is a Ricci-flat \emph{\aKa}metric
on $M$, compatible with the opposite orientation to the one induced on $M$ by
the \hk structure $(\wM_1,\wM_2,\wM_3)$. The integrability of the \acs $I$ is
equivalent to the flatness of the metric $\gM_0$, and this is the only
possibility when $M$ is compact, see \cite{SekSCE}. However, completely
explicit local examples of 4-dimensional \hk manifolds admitting a
\emph{non-integrable} \aKa structure $I$ are now known
\cite{ACGGWS,ArmAAK,NuPFDE}.\vs

\begin{thm}\label{thm2} Let $(M,\gM_0,\wM_1,\wM_2,\wM_3)$ be a \hk
4-manifold. Let $\k,\l$ be real constants with $\k+\l t>0$ for any $t\in\ab$
with $a>0$, and set \[\twM=(\k+\l t)\wM_1,\q\s=t(\k+\l t)^2.\] Then the 3-form
\textrm{(\ref{phi})} defines a torsion-free $G_2$ structure on $\ab\ti
\R_{\u,\v}^2\ti M'$ where $M'$ is a suitable open subset of $M$.

Suppose, furthermore, that $(M,\gM_0)$ admits an \aKa structure $(\wM_0,I)$
compatible with the opposite orientation to the one induced by
$(\wM_1,\wM_2,\wM_3)$. Let $p,q,\k,\l$ be real constants satisfying
\textrm{(\ref{inequality})} for $t\in\ab$, $a>0$, and set
\[\twM=(p+qt)\wM_0+ (\k+\l t)\wM_1,\q\s=t\Big((\k+\l t)^2-(p+qt)^2\Big).\]
Then \textrm{(\ref{phi})} again defines a torsion-free $G_2$ structure on a
manifold of the form $\ab\ti \R_{\u,\v}^2\ti M'$.\end{thm}\vs

\begin{rem} It suffices to take $M'$ to be any contractible open subset of
$M$, in which case the 1-forms $\xi_P$ and $\eta_N$ of Theorem~\ref{thm1} can
always be defined on $P=\R_t^+\ti M'$ and $N\=\R^+_t\ti\R_\u\ti M'$
respectively. However, as we shall see below, we can alternatively keep the
4-manifold $M$ fixed and think of the 1-forms $(\xi,\eta)$ in
Theorem~\ref{thm1} as connection forms of a principal $\T^2$ bundle $W$ over
$M$. Since \begin{equation}\label{dd}(d\xi,d\eta)=(q\wM_0+s\wM_1,\,-\wM_2)
\end{equation} is the curvature of the principal connection, we obtain {\it
integrality} constraints for the cohomological classes $[\fp(q\wM_0+\l\wM_1)]$
and $[\fp\wM_2]$ of $M$, in the sense that they must be contained in the image
of the universal morphism $H^2(M,\Z)\to H^2(M,\R)$.

Note also that in general only two of the four real parameters $(p,q,\k,\l)$
are effective in the sense that one can fix two of them, by rescaling $\varphi$
and $V$ on $Y$. Moreover, if $q=\l=0$, (\ref{phi}) corresponds to the product
of a Calabi-Yau 6-manifold with $\R_{\u}$, while generically the corresponding
metric has holonomy equal to $G_2$ in accordance with
Remark~\ref{irreducible}.\end{rem}

By way of proving Theorem~\ref{thm2}, we shall describe the reverse
construction in the case $(p,q,\k,\l)\=(0,0,0,1)$ for simplicity, so that
$\s=t^3$. This is also the case in which the family $g(t)$ consists of {\it
homothetic} hyperk\"ahler metrics on $M$. The resulting metrics with holonomy
$G_2$ were first described in \cite{GLPSDW} in the context of $\T^2$
bundles. To duplicate this situation, we assume that the 2-forms $\fp\wM_1$ and
$\fp\wM_2$ are integral. This is always true locally, though if $M$ is compact
(i.e.\ a K3 surface or a torus) the assumptions imply that $(M,J_3)$ is
\emph{exceptional} -- its Picard number is maximal, see e.g.\ \cite{BPVCCS}. On
any exceptional K3 surface, Yau's theorem implies the existence of \hk metrics
satisfying these integrality assumptions.

Let $P$ be the total space of the principal $S^1$ bundle over the \hk
4-manifold $M$ classified by $[\fp\wM_1]$, and $\xi$ a connection 1-form on $P$
such that $d\xi=\wM_1$ in accordance with (\ref{dd}). Attached to $P$ is also a
principal $\C^*$ bundle over $M$, whose total space $N$ is the real 6-manifold
manifold $N \cong \R_t^+\ti P$.  Pulling back forms to this new total space, we
define on $N$ an almost-Hermitian structure $(\gN,\wN,J)$ by \begin{equation}
\label{gN}\ba{c} \wN = t\wM_1+ dt \we\xi\y \gN =\ds t\kM+t^{-3}\xi\ot\xi+t^3dt
\ot dt\ea\end{equation} and \[J\al=J_1\al, \ \al\in\La^1M,\qq J\,dt=t^{-3}\xi.
\] The 2-form $\wN$ is closed, satisfies (\ref{gJ}), and it is easy to check
that $J$ is integrable. Thus, $(\gN,\wN,J)$ defines a \Ka structure on $N$.

Now let $U$ be vector field which is $\gN$-dual to $t^{-3}\xi$ so that $U$ is
tangent to the fibers of $P$ and $\xi(U)=1$; it follows that $U$ is the
generator of the natural $S^1$ action on $N$, acting as a rotation on each
fibre of $P$. Moreover, $U$ is a Hamiltonian isometry of the \Ka structure
(\ref{gN}) and the corresponding momentum map is $-t$.

The integrality assumption for $\wM_2$ implies that there exists a principal
$S^1$ bundle over $N$, classified by $[-\fp\wM_2]\in H^2(N,\Z)$. We denote by
$Y$ the corresponding 7-dimensional total space and take $\eta$ to be a
connection 1-form satisfying (\ref{deta}).\vs

\begin{cor}\label{Gibb} With the above assumptions, the $3$-form \[\varphi=
t\,\wM_1\we\eta+dt\we\xi\we\eta+t^4\wM_3\we dt+t\,\wM_2\we\xi\] defines a
torsion-free $G_2$ structure on $Y$. The corresponding 4-form is \[*\varphi=
\wM_3\we\xi\we\eta -t^3\wM_2 \we dt\we\eta+t^3\wM_1\we dt\we\xi+\ft12t^4
\wM_1\we\wM_1.\]\end{cor}\vs

\begin{rem} {\rm In the above situation, the \Ka structure $(\gN,\wN,J)$ on
$N$ belongs to the classes of metrics recently studied in \cite{ACGHFK} and
\cite{DeMLTC}. In the terminology of the former, $(\gN,J)$ arises from a
Hamiltonian 2-form of order 1 with one non-constant eigenvalue (equal to $t$)
and one zero constant eigenvalue of multiplicity 2; the structure function
$F(t)$ of $\gN$ is $t^{-1}$.

Moreover, \[\wN=t\wM_1+dt\we\xi=d(t\xi)=d(t^4 Jdt)= dd^c(\ft15t^5),\] using the
definition of $\xi$ after (\ref{deta}), and (\ref{dd}),(\ref{gN}). Not only
then is $f(t)=\ft12\log t$ a Ricci potential for $\gN$ (see
Corollary~\ref{cor2}), but $\ft25t^5=\ft25e^{5f}$ acts as a \Ka potential for
the same metric.}\end{rem}\vs

We now return to the general case of Theorem~\ref{thm2}, when $p$ and $q$ are
not both zero. Now the hyperk\"ahler metrics $g(t)$ are not homothetic, though
for the examples in \cite{ACGGWS,ArmAAK,NuPFDE} the family $g(t)$ is an isotopy
(i.e.\ $g(t)$ are all isometric to an initial hyperk\"ahler metric, under the
flow of a vector field on $M$).

If the cohomology classes of $\fp(q\wM_0+\l\wM_1)$ and $\fp\wM_2$
are integral, the $G_2$ structure is defined on the product of an
interval and the principal $\T^2$ bundle $W$ over $M$, classified
by $[\fp(q\wM_0+\l\wM_1)]$ and $[-\fp\wM_2]$. The action of the
2-torus $\T^2$ is generated by the commuting vector fields $U=
k(\cd,t^{-2}\xi)$ and $V=k(\cd,t^{-2}\eta)$, which preserve the
$G_2$ structure.\vs

Take $M=\T^4$ to be the 4-torus with a flat \hk metric. There are
three cases according to the signature of the bilinear form $B$ of
(\ref{B}) restricted to the 2-dimensional subspace generated by (\ref{dd}):
\begin{enumerate}

\item $\l>q$, so $B$ is positive definite and there is a flat metric $g_0$ on
$\T^4$ with $d\xi,d\eta\in\La^2_+$. Thus, there exists a basis $(e^i)$ of
1-forms on $\T^4$ such that $d\xi=\wM_1$ and $d\eta=-\wM_2$ with \[\wM_1=
e^{14}+e^{23}, \q\wM_2=e^{13}+e^{42},\q\wM_3=-(e^{12}+e^{34}),\] reflecting the
structure of the real Lie algebra $\mathfrak h$ associated to the complex
Heisenberg group $H$. Then $W=\Gamma\backslash H$ is the Iwasawa manifold
\cite{AGSAHG,AGSHGI}.

\item $\l=q$ so that $d\xi$ is a simple 2-form. In this case, we may assume
that $d\xi=e^{24}$ and either $d\eta=e^{14}+e^{23}$ or $d\eta=e^{14}$. Then $W$
is a $\T^2$ bundle over $\T^4$ corresponding to one of two other
nilmanifolds. The simplest holonomy $G_2$ metric in the first case is obtained
by setting $(p,q,\k,\l)=(-1,1,1,1)$ so that $\s=4t^2$.  The resulting metric
coincides with that described in Example 2 of \cite[\S4]{ChSITS} (except that
$t^2$ there has now become $t>0$).

\item $\l<q$ and there is a basis with $d\xi=e^{14}-e^{23}$ and $d\eta=
e^{14}+e^{23}$, so $\langle d\xi,d\eta\rangle=\langle d{\tilde\xi},d{\tilde
\eta}\rangle$ with $d{\tilde\xi}=e^{14}$ and $d{\tilde\eta}=e^{23}$. In this
case, we may take $W$ to be a discrete quotient of $H_3\ti H_3$, where $H_3$ is the real Heisenberg group.\end{enumerate}\vs

\begin{ex} To make (i) more explicit, we may take coordinates $\la,\mu,\ell,m$
on $\T^4$ and fibre coordinates $x,y$ on $W$ so that
\[e^5=-\eta=dx-\la d\ell+\mu dm,\q e^6=\xi=dy-\mu d\ell-\la dm\]
are the corresponding connection 1-forms. The resulting $G_2$ metric
\[ k=t^2(d\la^2+d\mu^2+d\ell^2+dm^2)+t^{-2}(dx-\la d\ell+\mu dm)^2
+t^{-2}(dy-\mu d\ell-\la dm)^2+t^4dt^2\] is defined on $(0,\infty)\ti W$.  It
is shown in \cite{GLPSDW} to arise from an $SO(5)$ invariant $G_2$ metric on
the total space of $\Lambda^2_+\to S^4$ by a contraction of the isometry group.

The explicit form above makes it easy to compute the the Riemann tensor
$R_{ijkl}$ of $k$. The package \texttt{GRTensor} (available from
\texttt{http://grtensor.phy.queensu.ca}) was in fact used to verify that (a)
the Ricci tensor $R^i{}_{jil}$ is zero, (b) the matrix $R^i{}_{jkl}$ with rows
labelled by $i$ has rank 7 everywhere, and (c) the matrix $R^{ij}{}_{kl}$ with
rows labelled by $(i,j)$ has rank 14 everywhere. Point (b) confirms that $k$ is
irreducible and therefore has holonomy \textit{equal} to $G_2$, which is
consistent with (c). The same technique can be used to analyse metrics arising
from (ii) and (iii).\end{ex}

\vs\section{Hypersurface structures}

We shall now investigate the metrics with holonomy $G_2$ constructed in
Theorem~\ref{thm2} by restricting them to hypersurfaces on which $t=
\|\eta\|_\gY$ is constant \textit{before} taking an $S^1$ quotient. These
hypersurfaces correspond to the total spaces $W$ of the $\T^2$ bundles of
Example~1.

We can write the 3-form (\ref{phi}) as \[\varphi=(\twM+dt\we\xi)\we\eta+
t(\s\,\wM_3\we dt+\wM_2\we\xi)\] where $\xi,\,\eta$ are the corresponding
connection 1-forms on $N,\,Y$ respectively. Since $\|dt\|^2_\gN=\s^{-1}$, it
follows from (\ref{gNgY}) that $\|dt\|^2_{\gY}=\x^{-1}$ where $\x=\s t$, so
\begin{equation}\label{x}\x=t^2\Big((\k+\l t)^2-(p+qt)^2\Big)\end{equation}
The 1-form $\x^{1/2}dt$ therefore has unit norm relative to the $G_2$ metric
$\gY$, and we may write \begin{equation}\label{tau}\x^{1/2}dt=d\tau
\end{equation} where \[\tau=\int\!\!t\Big((\k+\l t)^2-(p+qt)^2\Big)^{1/2}dt.\]
The important point here is that $\s$ is constant as a function on $M$, and so
$\x$ is really just a function of $t$.

As a consequence, \begin{equation}\label{standard} \ba{c}\varphi=\rho\we
d\tau+\phi^+\yy *\varphi =\phi^-\we d\tau+\ft12\rho\we\rho,\ea\end{equation}
where \begin{equation}\label{triple}
\ba{c}\rho=\x^{1/2}\wM_3+\x^{-1/2}\xi\we\eta,\yy
\phi^+=\twM\we\eta+t\wM_2\we\xi,\yy
\phi^-=t^{-1}\x^{1/2}\wM_2\we\eta-t^2\x^{-1/2}\twM\we\xi.\ea \end{equation}
Restricting to an interval where the function $t\mapsto\tau$ is bijective, let
$Y_\tau$ denote the hypersurface of $Y$ for which $\tau$ has the constant value
`$\tau$' (excusing the abuse of notation). Whereas $\wN,\psi^+,\psi^-$
characterize the $SU(3)$ structure on $N$, we are using $\rho,\phi^+,\phi^-$ (a
lexicographic shift) for the corresponding objects on $Y_\tau$.

The closure of the forms (\ref{standard}) is equivalent to asserting that
\begin{equation}\label{halfflat} d\phi^+=0,\qq d(\rho^2)=0,\end{equation}
($\rho^2$ denotes $\rho\we\rho$) for every fixed value of $\tau$, and that the
$SU(3)$ structures on $Y_\tau$ satisfy the equations
\begin{equation}\label{evolve}\frac{\pd\phi^+}{\pd\tau}=d\rho,\qq
\frac{\pd}{\pd\tau}(\ft12\rho^2)=-d\phi^-.\end{equation} An $SU(3)$ structure
satisfying (\ref{halfflat}) is called \textit{half-integrable} or
\textit{half-flat} in the formalism of \cite{ChSITS}.

To verify (\ref{evolve}) directly in our situation, first observe that
\[d(\xi\we\eta)=(q\wM_0+\l\wM_1)\we\eta+\xi\we\wM_2\] (recall (\ref{dd})),
and \begin{equation}\label{pair}\ba{l}\phi^+=(p\wM_0+\k\wM_1)\we\eta+t\,
d(\xi\we\eta)\yy\ft12\rho^2=\xi\we\eta\we\wM_3+\ft12\x\,\wM_3\we\wM_3.\ea
\end{equation} These equations explain the significance of the coordinates
$(t,\x)$ -- the above forms are linear in $t$ and $\x$ respectively. The reader
may now check that (\ref{evolve}) implies both (\ref{x}) and (\ref{tau}); a
constant of integration may be absorbed into the term $\k^2-p^2$.

Hitchin discovered that (\ref{evolve}) leads to a Hamiltonian system in the
symplectic vector space ${\bf V}\ti {\bf V}^*$ where ${\bf V}$ is the space of exact 3-forms on
the compact manifold $Y_{\tau}$, whose dual ${\bf V}^*$ can be identified with the
space of exact 4-forms. The Hamiltonian function $H$ is derived from
integrating volume forms determined algebraically by $\phi^+$ and $\rho^2$,
and this also enables $\phi^-$ to be determined from $\phi^+$ in
(\ref{triple}).  Elements of ${\bf V}$ represent deformations of $\phi^+$ in a fixed
cohomology class, and those of ${\bf V}^*$ deformations of $\rho^2$. Given a
solution of (\ref{halfflat}) for $\tau=a$, a solution of (\ref{evolve}) can
then be found on some interval $\ab$ \cite{HitSFS}. This approach also
underlies some of the newly constructed metrics (such as \cite{BGGGTL}) with
reduced holonomy.

The function $H$ is already implicit in our calculations above, which may be
summarized in the following way:

\begin{prop}\label{Ham} The solution (\ref{triple}) can be expressed in the
form $H=0$ where the function
$H=2t\Big((\k+\l t)^2-(p+qt)^2\Big)^{1/2}\! -2{\x}^{1/2}$ satisfies
\[\frac{dt}{d\tau}=-\frac{\pd H}{\pd {\x}},\qq \frac{d \x}{d\tau}=\frac{\pd
H}{\pd t}.\]\end{prop}\vs

It is an important consequence of the \Ka assumption that, for each choice
$(p,q,\k,\l)\in\R^4$, there is only one valid solution curve in the $(t,\x)$
plane.\vs

\begin{ex} In the light of Proposition~\ref{Ham}, a variant of (\ref{triple})
and Example~1 is provided by setting
\[\ba{c}\rho=\pm\x^{1/2}(e^{12}+e^{34})+\x^{-1/2}e^{56},\y
\phi^+=\phi^+_0+t\,d(e^{56}),\y\phi^-=\phi^-_0\pm t(e^5\we de^5+e^6\we de^6),
\ea\] with \begin{equation}\label{quartic} \x=(t^2-\ft12H)^2.\end{equation} We
assume that $\x^{1/2}>0$, so that the orientation $\rho^3$ remains fixed. The
notation is consistent with that of \cite{AGSAHG}, with $de^5=e^{13}+e^{42}$
and $de^6=e^{14}+e^{23}$. Any value of the constant $H$ gives a valid solution
and so a holonomy $G_2$ metric on the product $\ab\ti W$ of some interval
with the Iwasawa manifold $W$, though $\phi^\pm_0$ and the signs must be chosen
to ensure compatibility (essentially (\ref{comp2})) and positive definiteness
of the resulting metric.

The Figure plots the quartic curves (\ref{quartic}) in the $(t,\x)$ plane for
various values of $H$; two pass through each point because of the sign
ambiguity implicit in the definition of $H$.

\begin{enumerate}\item $H=2$ gives $\x=(1-t^2)^2$, including the bell-shaped
segment. To satisfy (\ref{evolve}) for $|t|<1$ we need to choose the plus sign
in $\phi^-$. We may then define the 3-forms $\phi^\pm_0$ by setting
\[\phi^++i\phi^-=i((e^5+ie^6)+t(e^5-ie^6))\we d(e^5+ie^6).\] It follows that
$\phi^+_0+i\phi^-_0$ is a form of type $(3,0)$ relative to the standard complex
structure on $W$, and taking $\x=1$ and the plus sign in $\rho$ determines a
compatible Hermitian metric via (\ref{gJ}). As one flows away from the point
$(0,1)$, the \acs induced on the hypersurface becomes non-integrable;
$\rho,\phi^\pm$ degenerate simultaneously when one reaches
$t=\pm1$. Furthermore, we may take $\tau=\ft13t^3+t$, a cubic equation with
solution \[t=\al^{1/3}-\al^{-1/3},\q 2\al= 3\tau+(9\tau^2-4)^{1/2},\] in
contrast with the \Ka scheme in which $\tau=c+\ft12\k t^2+\ft13\l t^3$ in the
case $p=q=0$.

\FFF

\item $H=-2$ gives $\x=(t^2+1)^2$ (the curve above and touching the bell). This
requires the minus sign in $\phi^-$ and provides a different deformation of the
standard Hermitian structure on $W$. Indeed, $\phi_0^++i\phi_0^-$ is modified
by the addition of a form of type $(1,2)$ rather than $(2,1)$, and the
resulting \acs is undefined when $|t|$ reaches 1. However, $\rho$ remains
non-degenerate for all $t$, and this corresponds to a different singular
behaviour of the $G_2$ metric.

\item $H=0$ and $\x=t^4$ (with the flattened base) requires $\phi^\pm_0=0$ and
the minus signs in $\phi^-,\rho$, reproducing exactly the solution of
Example~1. The \acs induced on $W$ by $\phi^+$ is \textit{constant} and was
first singled out for study in \cite{AGSHGI}, where it is called $J_3$. The
2-form $\rho$ degenerates only for $t=0$, and the resulting metric has the
advantage of being `half-complete'.
\end{enumerate}\end{ex}\vs

We conclude this section by providing an explanation of why the system
(\ref{evolve}) often reduces to one variable. Let $\phi^+$ be an invariant
closed 3-form with stabilizer $SL(3,\C)$ on either of the four nilmanifolds
$\Gamma\backslash H$ described just before Example~1. Let $\phi^-$ be the
unique 3-form for which $\phi^+\+i\phi^-$ has type (3,0) relative to the \acs
$J$ determined by $\phi^+$. Let $\ker d$ denote the space of invariant closed
1-forms (equivalently, the kernel of the natural mapping $d:\mathfrak{h}^*\to
\alt^2\mathfrak{h}^*$).\vs

\begin{lemma} In the above situation, $d\psi_+=0$ implies that $J(\ker d)=
\ker d$. \end{lemma}\vs

\noindent This can be proved by generalizing an argument in the proof of
\cite[Theorem~1.1]{KeSCSI}, which draws the same conclusion if $J$ is
integrable. In view of the structure equations for the Lie algebra $\mathfrak
h$, the condition that $J$ leave $\ker d$ invariant is equivalent to asserting
that $d\phi^-$ belong to the 1-dimensional space \[\langle e^{1234}\rangle=
\alt^4(\ker d),\] but it is this fact that simplifies the equations.

\vs\section{Solutions varying on $M$}

As another ramification of the evolution equations for $(\twM\=\twM(t),\,\s\=
\s(t))$, we consider the case when $(M,\gM_0,\wM_1,\wM_2,\wM_3)$ is a \hk
manifold, and \[\ba{c}\twM=\wM_1-\ft12d_M d_M^c\,G,\yy 2\s=G''\ea\] for a
smooth function $G$ on $\ab\ti M$, where $\>'\>$ continues to denote $\pd/\pd
t$. We assume here that $\twM$ is a positive definite (1,1)-form on $(M,J_1)$
for each fixed $t$ in the interval $\ab$. Whilst (\ref{flow}) is
automatically satisfied, (\ref{Tcompatible}) becomes \begin{equation}
\label{monge-ampere} 2t\,\cM(G)= G''\end{equation} where $\cM$ denotes the
complex Monge-Amp\`ere operator on $(M,\gM_0,J_1)$, defined by \[\Big(\wM_1-
\ft12d_Md_M^cf\Big)^{\we 2}=\cM(f)\wM_1\we\wM_1,\] for all $f\in C^\infty(M)$.
Note that the 1-form $\xi_P$ of Theorem \ref{thm1} is automatically defined on
$P=\R_t^+\ti M$ by \begin{equation}\label{xi-t}\xi_P=-\ft12 d_M^c(G').
\end{equation}

Since the \hk metric $\gM_0$ is Ricci-flat and the $\wM_i$'s are parallel
2-forms, there exists a real-analytic structure on $M$, compatible with
$(\gM_0,J_1,\wM_1,\wM_2,\wM_3)$. Thus, applying the Cauchy-Kowalewski theorem,
one obtains a general existence result.\vs

\begin{cor}\label{real-an} Let $(M,\gM_0,\wM_1,\wM_2,\wM_3)$ be a
hyperk\"ahler, real-analytic 4-manifold and $J_1$ be the \Ka structure
compatible with $\wM_1$. Suppose that $G^0,G^1$ are real-analytic functions on
$M$ with $\wM_1-\ft12 d_M d^c_M G^0$ positive definite with respect to
$J_1$. Then there exist a real number $a>0$ and a real-analytic solution
$G(t,\cd)$ of (\ref{monge-ampere}) defined on $(0,a)\ti M$ with $G(0)=G^0$,
$G'(0)=G^1$, and such that $\wM_1-\ft12 d_Md^c_M G$ is positive definite for
any $t\in(0,a)$. Thus, $\twM=\wM_1-\ft12 d_M d_M^c G$ and $\s=\ft12G''$
define, via Theorem \ref{thm1}, a torsion-free $G_2$ structure on a manifold
of the form $(0,a)\ti \R^2_{\u,\v}\ti M'$ where $M'$ is a suitable open subset
of $M$.\end{cor}\vs

In the above construction $M'$ should be taken so as to solve $d\eta_N= -\wM_2$
for a 1-form $\eta_N$ on $N=\R^+_t\ti\R_{\u}\ti M'$ (see Theorem~{\ref{thm1}
and (\ref{xi-t})). Alternatively, we may assume that the cohomology class
$[\fp\wM_2]$ of $M$ is integral and $\eta$ is a principal connection of the
principal $S^1$ bundle $Q$ over $M$, classified by $[-\fp\wM_2]$. In this case
Corollary~\ref{real-an} produces examples of torsion-free $G_2$ structures with
a ${\mathbb R}\times S^1$ symmetry on $Y=(0,a)\ti\R_{\u}\ti Q$. One has no control over
the real number $a$.

It is tempting to spot some special solutions to (\ref{monge-ampere}), by
reducing the problem to a linear (elliptic) equation. This can be done by
assuming that for each $t>0$ the function $G$ generates a complex {\it
Monge-Amp\`ere foliation} on $(M,J_1)$ (see e.g.~\cite{BeKa}), meaning that
$$d_M d^c_M G\we d_M d^c_M G=0$$ and $d_M d^c_M G$ has constant rank. (The
integral curves of $d_M d^c_M G$ then foliate $M$ by complex submanifolds.)
The point is that in this case we have \[\cM(G)=1+\ft12\Delta G,\] where
$\Delta$ is the Riemannian Laplacian of the \hk metric $\gM_0$.

The above situation appears in particular when $(M,J_1)$ admits a holomorphic
$\C$ action and we look for equivariant solutions of (\ref{monge-ampere}), or
when $(M,J_1)$ admits a holomorphic fibration $p:M \to C$ over a complex curve
$C$ and we look for solutions of the form $G\circ p$ where for each $t$, $G$
is a function on $C$.

Consider finally the flat \hk metric $\gM_0$ on $(M,J_1)=\C^2\cong\R^4$,
determined by the 2-forms \[\wM_1=\ft12i(dz_1\we d\ol z_1+dz_2\we d\ol z_2),
\q\wM_2 =\Re(dz_1\we dz_2),\q\wM_3=\Im(dz_1\we dz_2),\] where $z_1=\la+i\mu$
and $z_2=\ell+im$ are the canonical coordinates of $\C^2$. Letting $G(t,z_1)$
be a function on $\C^2$ which does not depend on $z_2$, the equation
(\ref{monge-ampere}) reduces to \[G''+ t\left(\frac{\pd^2 G}{\pd\la^2}+
\frac{\pd^2 G}{\pd\mu^2}\right)=2t.\] Separable solutions are given by
$G(t,\la,\mu)=\ft13t^3+H(\la,\mu)K(t)$, where \[\frac{\pd^2 H}{\pd\la^2}+
\frac{\pd^2H}{\pd\mu^2}+cH=0,\] $c$ is a constant and $K$ a solution of the
Airy equation $K''=ctK$ (equivalently, $L=K^{-1}K'$ satisfies the Riccati
equation $L'+L^2=ct$).

\begin{ex} Taking $H$ to be periodic in $\la,\mu$ (which requires $c>0$)
yields solutions of (\ref{monge-ampere}) defined on $(0,a)\times\T^4$. One such
example is obtained by taking $H=\sin\la$ (so $c=1$), and \[K=\Ai(t)=\ft13
t^{1/2}(J_{1/3}(\zeta)+J_{-1/3}(\zeta)),\q \zeta=\ft23it^{3/2}.\] Setting
$f=1+\Ai(t)\sin\la$, the resulting $G_2$ metric \[t(fd\la^2+fd\mu^2+d\ell^2
+dm^2)+f^{-1}(dx-\Ai'(t)\cos\la\,d\mu)^2 +t^{-2}(dy-\la d\ell+\mu
dm)^2+t^2f\,dt^2\] is Ricci-flat and irreducible. Since $\Ai(t)\to0$ as
$t\to\pm\infty$, the above metric is asymptotic to a constant solution with
$u=t$ and holonomy equal to $SU(3)$ (see Remark~3). However, the above
construction can be easily modified to provide explicit deformations of the
non-trivial metrics constructed in \S3.\end{ex}

\enddocument